\documentclass[12pt]{article}

\usepackage{graphicx}
\usepackage{amsmath}

\newcommand{\eps}{\varepsilon}
\def\qed{\hfill $\squar$}
\def\squar{\vbox{\hrule\hbox{\vrule height 6pt \hskip 6pt\vrule}\hrule}}

\usepackage{mathtext}
\usepackage[T2A]{fontenc}
\usepackage[cp866]{inputenc}

\usepackage[english,russian]{babel}
\usepackage{amsfonts,amssymb}
\usepackage{latexsym}
\usepackage{personal}
\newtheorem{thm}{Теорема}
\newtheorem{lem}{Лемма}

\newtheorem{pro}{Предложение}

\newtheorem{opr}{Определение}

\begin{document}

\renewcommand{\figurename}{Фиг.}
\large

УДК 517.928.7; 531.8; 531.36

\setcounter{equation}{0} \thispagestyle{empty}
\begin{center}
{\bf Асимптотическая устойчивость колебаний двухмассного
резонансного грохота с односторонним ограничителем без зазора\\}

О. Ю. Макаренков
\end{center}

{\bf Аннотация.} В работе доказывается асимптотическая
устойчивость периодических колебаний в модели резонансного грохота
в предположении, что частоты порождающей системы соотносятся как
1:2 и частота внешнего двигателя совпадает с меньшей из этих
частот. Такая постановка соответствует широко используемому режиму
работы грохота -- нелинейному резонансу. Обоснование проводится
при помощи предложенной автором негладкой версии второй теоремы
Н.Н. Боголюбова. Строго доказано, что найденный резонансный режим
является двухчастотным.

{\bf 1. Введение.} В работе изучается поведение двухмассного
резонансного грохота, расчетная схема которого представлена на
фиг.~1.

\setcounter{subsection}{1}\setcounter{equation}{0}

\begin{figure}[h]
\begin{center}
  \includegraphics[scale=.5]{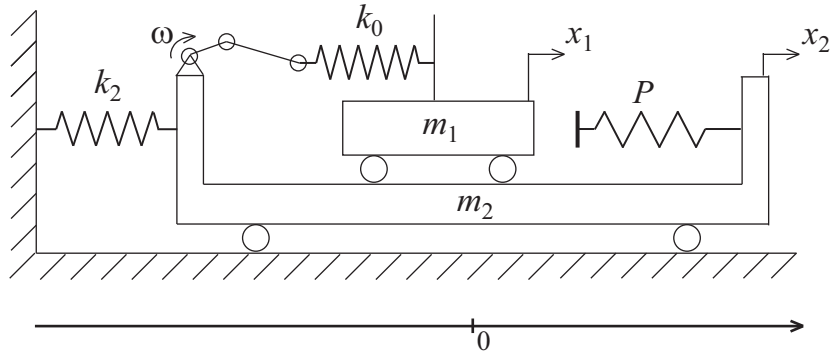}\\
  \caption{\footnotesize Схема двухмассного резонансного грохота, приводимого в движение
двигателем с частотой вращения $\omega$.  Жесткость пружины $P$
описывается кусочно линейной функцией
$P(x_1-x_2)=\max\left\{0,x_1-x_2\right\}$. В состоянии покоя и при
выключенном двигателе зазор между ограничителем $P$ и телом $m_1$
равен нулю.}
  \end{center}
\end{figure}

\vskip-0.7cm

При выключенном двигателе и отсутствии трения система совершает
двухчастотные колебания с частотами $\omega_1$ и $\omega_2$. Если
включить двигатель и считать трение малым, то приближая $\omega$ к
$\omega_1$ или $\omega_2$, в системе возникает резонанс. В
\cite[с.~137]{kru} замечено, что один из этих резонансов
$w=\omega_2$ является нелинейным и что причиной нелинейности
является содержащийся в модели упругий ограничитель. Для
обоснования данного утверждения в книге \cite{kru} использован так
называемый метод эквивалентной линеаризации, являющийся не вполне
строгим.

В настоящей работе показывается, что при формальном требовании
малости жесткости ограничителя указанный нелинейный резонанс может
быть строго изучен методом усреднения и, что асимптотическая
устойчивость соответствующих колебаний может быть установлена при
помощи негладкого аналога второй теоремы Н.Н. Боголюбова
\cite{blm}.

Статья построена следующим образом. В разделе 2 разрабатывается
способ приведения к стандартной форме принципа усреднения системы
дифференциальных уравнений, описывающих наиболее общую
механическую систему двух связанных тел. В разделе 3 вводятся
условия применимости негладкого аналога второй теоремы Н.Н.
Боголюбова. В третьем разделе полученный способ применяется к
системе дифференциальных уравнений, описывающих модель грохота
фиг.~1, в которой для простоты предполагается
$\omega_2=2\omega_1=2\omega$. Доказывается, что наличие упругого
ограничителя приводит к асимптотически устойчивым периодическим
колебаниям с двумя частотами $\omega$ и $2\omega$.

{\bf 2. Общее уравнение для амплитуды асимптотически устойчивых
периодических колебаний в слабонелинейных механических системах,
содержащих два тела.} \setcounter{subsection}{2}
\setcounter{equation}{0} Систему обыкновенных дифференциальных
уравнений, описывающих координаты $x_1$ и $x_2$ двух тел,
связанных посредством слабо нелинейных пружин, при учете вязкого
трения и периодического воздействия, в наиболее общем виде можно
записать как
\begin{equation}\label{ps1}
\begin{pmatrix}
  m_1 \ddot{x}_1\\m_2\ddot{x}_2
\end{pmatrix}  +Mx+
\begin{pmatrix}
 \varepsilon \dot{q}_1(Q_1(\varepsilon)x)Q_1(\varepsilon)\dot{x}\\\varepsilon \dot{q}_2(Q_2(\varepsilon)x)Q_2(\varepsilon)\dot{x}
\end{pmatrix}=\varepsilon F(t,x,\dot{x},\varepsilon)
\end{equation}
где $m_1, m_2 > 0$ -- массы тел, $\varepsilon>0$ -- малый
параметр, $M$ -- $2 \times 2$-матрица, $Q_1(\varepsilon),
Q_2(\varepsilon) - 1 \times 2$-матрицы и $F: \mathbb{R} \times
\mathbb{R}^2 \times \mathbb{R}^2 \times [0,1] \to \mathbb{R}^2$.
Условия на использованные функции будут возложены позже, отметим
лишь, что скалярные функции $q_1$ и $q_2$  непрерывны и
кусочно-линейны, так что соответствующие производные $\dot{q}_1$ и
$\dot{q}_2$ -- кусочно-постоянны. Совершим ряд преобразований.
Перепишем систему (\ref{ps1}) как
\begin{equation*}
\begin{pmatrix}
  m_1 \ddot{x}_1\\m_2\ddot{x}_2
\end{pmatrix}  +Mx+ \varepsilon
\begin{pmatrix}
{(q_1(Q_1(\varepsilon) x(\cdot)))}^{'} \\ {(q_2(Q_2(\varepsilon)
x(\cdot)))}^{'}
\end{pmatrix}=\varepsilon F(t,x,\dot{x},\varepsilon)
\end{equation*}
откуда видно, что вводя новые переменные (см. \cite{lev})
\begin{equation*}
\begin{matrix}
  y_1(t)=m_1 \dot{x}_1(t)+\varepsilon q_1(Q_1(\varepsilon)x(t)),\
  \ \
  y_2(t)=m_2 \dot{x}_2(t)+\varepsilon q_2(Q_2(\varepsilon)x(t))
\end{matrix}
\end{equation*}
можем переписать ее в виде четырех уравнений первого порядка
\begin{eqnarray}
  \dot{x}_1&=&\frac{1}{m_1}(y_1-\varepsilon
  q_1(Q_1(\varepsilon)x)),\ \ \
  \dot{x}_2=\frac{1}{m_2}(y_2-\varepsilon q_2(Q_2(\varepsilon)x)) \nonumber \\
  \begin{pmatrix} \dot{y}_1\\\dot{y}_2 \end{pmatrix}&+&Mx= \varepsilon
  F \left ( t,x,  \begin{pmatrix} \frac{1}{m_1}(y_1-\varepsilon q_1(Q_1(\varepsilon)x)) \\ \frac{1}{m_2}(y_2-\varepsilon q_2(Q_2(\varepsilon)x)) \end{pmatrix},\varepsilon \right ) \label{ps2}
\end{eqnarray}
Соответствующая порождающая система выписывается как
\begin{eqnarray}
  &&\dot{x}_1=\frac{1}{m_1}y_1, \qquad \dot{x}_2=\frac{1}{m_2}y_2 \label{np2} \\
  &&\dot{y_1}=-a_{11}x_1-a_{12}x_2,\qquad
  \dot{y_2} =-a_{21}x_1-a_{22}x_2 \label{np4}
\end{eqnarray}
где $a_{ij}$ -- компоненты матрицы $M.$ Для применения принципа
усреднения требуется предположить, что все решения системы
(\ref{np2})-(\ref{np4}) являются периодическими. Выясним условия,
при которых это имеет место. Подставляя первое уравнение
(\ref{np2}) в первое уравнение (\ref{np4}), второе уравнение
(\ref{np2}) во второе уравнение (\ref{np4}), и полученное в первом
случае уравнение в полученное во втором случае, получаем следующее
уравнение для $x_1$
\begin{eqnarray}\label{npx1}
  & &\ddddot{x}_1 +\frac{a_{11}}{m_1}\ddot{x}_1+\frac{a_{22}}{m_2}\ddot{x}_1-\frac{a_{21}a_{12}}{m_1m_2}x_1+\frac{a_{22}a_{11}}{m_1m_2}x_1=0
\end{eqnarray}
при этом из первого уравнения (\ref{np4})
\begin{equation}\label{npx2}
  x_2=\left({-a_{11}x_1-m_1\ddot{x}_1}\right)/{a_{12}}
\end{equation}
Поэтому, для того, чтобы система (\ref{np2})-(\ref{np4}) имела
только периодические решения необходимо и достаточно, чтобы
собственные значения характеристического многочлена уравнения
(\ref{npx1}) были различны и чисто мнимые. Данный
характеристический многочлен выписывается как \centerline{$
  \lambda^4+(a_{22}({1}/{m_2})+a_{11}({1}/{m_1}))\lambda^2 + (a_{11}a_{22}-a_{21}a_{12})/({m_1m_2})=0
$} и его корни которого находятся по формуле
\begin{equation}
  \lambda^2_{1,2}=
  -\frac{1}{2}\Big(a_{22}\frac{1}{m_2}+a_{11}\frac{1}{m_1}\Big)\pm \frac{1}{2} \sqrt{\Big(a_{22}\frac{1}{m_2}-a_{11}\frac{1}{m_1}\Big)^2+4\frac{1}{m_1m_2}a_{21}a_{12}} \label{LAM}
\end{equation}
Поэтому, каждое решение системы (\ref{np2})-(\ref{np4}) является
периодическим тогда и только тогда, когда
\begin{equation}\label{A1}
   a_{21}a_{12}<a_{11}a_{22}
\end{equation}
При этом, согласно правилу составления решений линейных уравнений
(см. \cite{pont}), общее решение уравнения (\ref{npx1})
записывается в виде $
  x_1(t)=a_{12}\left(A_{1C}\sin \omega_1t+A_{1S}\cos \omega_1t+A_{2C}\sin \omega_2 t+A_{2S}\cos
  \omega_2t\right),
$ где $w_{1,2}=\sqrt{-(\lambda_{1,2})^2}$ и
 $A_{1C}, A_{1S}, A_{2C}, A_{2S}$ -- произвольные постоянные.
Тогда, из (\ref{npx2}) имеем $x_2(t)
=(-a_{11}+m_1\omega^2_1)(A_{1C}\sin \omega_1t + A_{1S}\cos
\omega_1t)+(-a_{11}+m_1\omega^2_2)(A_{2C}\sin \omega_2t+A_{2S}\cos
\omega_2t).$ Таким образом, общее решение $(x_1, x_2, y_1, y_2)$
системы (\ref{np2})-(\ref{np4}) имеет вид: \centerline{$
  (x_1(t),x_2(t),y_1(t),y_2(t))^*= \Omega(t) (A_{1C},A_{1S},A_{2C},A_{2S})^*
$} где
\begin{equation*}
  \hskip-0.7cm \Omega(t)= \left( \begin{array}{ll}  a_{12}\sin(\omega_1t)&a_{12}\cos(\omega_1t)\\
  (-a_{11}+m_1\omega^2_1)\sin(\omega_1t)&(-a_{11}+m_1\omega^2_1)\cos(\omega_1t)\\
  m_1\omega_1a_{12}\cos(\omega_1t)&-m_1\omega_1a_{12}\sin(\omega_1t)\\
  m_2\omega_1(-a_{11}+m_1\omega^2_1)\cos(\omega_1t)&-m_2\omega_1(-a_{11}+m_1\omega^2_1)\sin(\omega_1t) \end{array} \right.
\end{equation*}
\begin{equation}\label{Omega}
\left. \begin{array}{ll}
a_{12}\sin(\omega_2t)&a_{12}\cos(\omega_2t)\\
(-a_{11}+m_1\omega^2_2)\sin(\omega_2t)&(-a_{11}+m_1\omega^2_2)\cos(\omega_2t)\\
m_1\omega_2a_{12}\cos(\omega_2t)&-m_1\omega_2a_{12}\sin(\omega_2t)\\
m_2\omega_2(-a_{11}+m_1\omega^2_2)\cos(\omega_2t)&-m_2\omega_2(-a_{11}+m_1\omega^2_2)\sin(\omega_2t)
\end{array} \right)
\end{equation}
\begin{lem}\label{lem1}{\bf \hskip-0.2cm.} Пусть $\Omega(t)$- фундаментальная система решений линейной системы
$  \dot{v}=Cv, $ где $C$ -- $n \times n$-матрица. Непрерывно
дифференцируемая функция $t \mapsto v(t)$  является на отрезке
$[0,t_0]$ решением системы
\begin{equation}\label{RSI}
  \dot{v}=Cv+g(t,v)
\end{equation}
где $g \in C^0(\mathbb{R}\times \mathbb{R}^n,\mathbb{R}^n)$, тогда
и только тогда, когда функция
\begin{equation}\label{ZAZA}
  u(t)=\Omega(t)^{-1}v(t)
\end{equation}
является на этом отрезке решением системы
\begin{equation}\label{RSISI}
  \dot{u}=\Omega(t)^{-1}g(t,\Omega(t)u)
\end{equation}
\end{lem}

Доказательство проводится аналогично \cite[лемма 2.4]{mak}.

Таким образом, для приведения системы (\ref{ps2}) к стандартной
форме принципа усреднения необходимо вычислить обратную к
(\ref{Omega}) матрицу. Прямой подстановкой проверяется, что эта
матрица дается формулой
\begin{equation*}
\Omega(t)^{-1}=
\left( \begin{array}{cccc} -\frac{1}{\omega_1}\cos \omega_1t&-\frac{1}{\omega_2}\sin \omega_1t&0&0\\
                                    \frac{1}{\omega_1}\sin \omega_1t&-\frac{1}{\omega_2}\cos \omega_1t&0&0\\
                                    0&0&\frac{1}{\omega_2}\cos \omega_2t&\frac{1}{\omega_1}\sin \omega_2t\\
                                    0&0&-\frac{1}{\omega_2}\sin \omega_2t&\frac{1}{\omega_1}\cos \omega_2t \end{array} \right) \circ
\end{equation*}
\begin{equation*}
\circ \frac{1}{(m_1)^2m_2a_{12}((\omega_1)^2-(\omega_2)^2)}\left( \begin{array}{cccc} \cos \omega_2t&-\sin \omega_2t&0&0\\
                                    \sin \omega_2t&\cos \omega_2t&0&0\\
                                    0&0&\cos \omega_1t&-\sin \omega_1t\\
                                    0&0&\sin \omega_1t&\cos \omega_1t \end{array} \right) \circ
\end{equation*}
\begin{equation*}
\circ \hskip-0.2cm\left(\hskip-0.2cm \begin{array}{cccc} m_2p_{22}m_1\omega_2\sin\omega_2t&-m_1a_{12}m_2\omega_2\sin\omega_2t&m_2p_{22}\cos\omega_2t&-m_1a_{12}\cos\omega_2t\\
m_2p_{22}m_1\omega_2\cos\omega_2t&-m_1a_{12}m_2\omega_2\cos\omega_2t&-m_2p_{22}\sin\omega_2t&m_1a_{12}\sin\omega_2t\\
m_2p_{21}m_1\omega_1\sin\omega_1t&-m_1a_{12}m_2\omega_1\sin\omega_1t&m_2p_{21}\cos\omega_1t&-m_1a_{12}\cos\omega_1t\\
m_2p_{21}m_1\omega_1\cos\omega_1t&-m_1a_{12}m_2\omega_1\cos\omega_1t&-m_2p_{21}\sin\omega_1t&m_1a_{12}\sin\omega_1t
\end{array} \hskip-0.2cm\right)\hskip-0.2cm
\end{equation*}
где $p_{21}=-a_11 m_1 \omega^2,$ $p_{22}=-a_{11}+m_1 \omega_2^2.$
Используя лемму~\ref{lem1}, можем записать систему (\ref{ps2}) в
стандартной форме принципа усреднения
\begin{equation}\label{AVE}
  \dot{A}=\varepsilon\Omega(t)^{-1}g(t,A,\varepsilon)
\end{equation}
где $A=(A_{1C}, A_{1S}, A_{2C}, A_{2S}) \in \mathbb{R}^4$,
\begin{equation}\label{g}
g(t,A,\varepsilon)=
\begin{pmatrix}
  -\frac{1}{m_1}q_1\left( Q_1(\varepsilon) \begin{pmatrix} \Omega_1(t), \Omega_2(t) \end{pmatrix}^* A \right)\\
  -\frac{1}{m_2}q_2\left( Q_2(\varepsilon) \begin{pmatrix} \Omega_1(t), \Omega_2(t) \end{pmatrix}^* A \right)\\
  \hskip-0.1cm F \hskip-0.1cm\begin{pmatrix}\hskip-0.1cm     t,\hskip-0.1cm \begin{pmatrix}\hskip-0.1cm
   \Omega_1(t)\hskip-0.1cm\\ \hskip-0.1cm\Omega_2(t)\hskip-0.1cm \end{pmatrix}\hskip-0.1cm A,\hskip-0.1cm
  \begin{pmatrix}
   \hskip-0.1cm\frac{1}{m_1}\hskip-0.1cm \left(\hskip-0.1cm\Omega_3(t)A- \varepsilon q_1\hskip-0.1cm \left(\hskip-0.1cm Q_1(\varepsilon)\hskip-0.1cm \begin{pmatrix}\hskip-0.1cm \Omega_1(t)\\ \Omega_1(t) \end{pmatrix}\hskip-0.1cm A\hskip-0.1cm \right)\hskip-0.3cm \right) \\
   \hskip-0.1cm\frac{1}{m_2}\hskip-0.1cm \left(\hskip-0.1cm\Omega_4(t)A- \varepsilon q_2\hskip-0.1cm \left(\hskip-0.1cm Q_2(\varepsilon)\hskip-0.1cm \begin{pmatrix}\hskip-0.1cm \Omega_1(t)\\ \Omega_1(t) \end{pmatrix}\hskip-0.1cm A\hskip-0.1cm \right)\hskip-0.3cm \right)
  \hskip-0.3cm\end{pmatrix}\hskip-0.3cm \end{pmatrix}\hskip-0.3cm
\end{pmatrix}
\end{equation}
и $\Omega_i(t)$ обозначает $i$-ю компоненту 4-мерного вектора
$\Omega_i(t)A$. Так как матрица $\Omega(t)$ равномерно ограничена
на $\mathbb{R}$, то асимптотически устойчивому решению $t \to
A(t)$ системы (\ref{AVE}) соответствует асимптотически устойчивое
решение $(x_1(t), x_2(t), y_1(t), y_2(t))^*=\Omega(t)A(t)$ системы
(\ref{ps2}). Более того, компоненты $(y_1, y_2)$ этого решения,
очевидно, дифференцируемы по времени и, значит, компонента $(x_1,
x_2)$ является, в этом случае, асимптотически устойчивым решением
исходной системы $(\ref{ps1})$.

Итак, установлено, что если выполнено условие (\ref{A1}), то
асимптотически устойчивым решениям $A(t)$ системы (\ref{AVE})
соответствуют асимптотически устойчивые решения
$(x_1(t),x_2(t))=(\Omega_1(t)A(t),\Omega_2(t)A(t))$ системы
(\ref{ps1}), где матрица $\Omega(t)$ дается формулой
(\ref{Omega}).

{\bf 3. Вторая теорема Н.Н. Боголюбова для систем с интегрально
дифференцируемыми правыми частями.} В этом разделе доказывается,
что правая часть полученной в разделе 2 системы (\ref{AVE})
удовлетворяет требованиям предложенного автором в \cite{blm}
негладкого аналога второй теоремы Н.~Н.~Боголюбова.
\setcounter{subsection}{3}\setcounter{equation}{0}
\begin{opr}\label{opr1}{\bf \hskip-0.2cm.}
Функцию $f \in C^0([0,T]\times \mathbb{R}^n \times
[0,1],\mathbb{R}^n)$ назовем интегрально дифференцируемой, если
каждому $\gamma
>0$ соответствует $\delta
>0$ и множество $M \subset [0,T]$, лебегова мера которого не
превосходит $\gamma
>0$, такие, что для всех $\| v-v_0\|<\delta, \quad t \in
[0,T]\setminus M$ и $\varepsilon \in [0,\delta]$ функция
$f(t,\cdot,\varepsilon)$ дифференцируема в точке $v$  и
\begin{equation*}
  \|f'_v(t,v,\varepsilon)-f'_v(t,v_0,0)\| \leqslant \gamma
\end{equation*}
\end{opr}

\begin{lem}\label{lem2}{\bf \hskip-0.2cm.}
Пусть $h \in C^0(\mathbb{R}\times \mathbb{R}^n \times
[0,1],\mathbb{R}^n)$ и\\ \centerline{$R^n_s=\{\xi \in
\mathbb{R}^n: {\rm sign} (\xi_i)=s_i, i \in \overline{1,n}\}, s
\in\{-1,1\}^n=\{-1,1\} \times \dots \times \{-1,1\}. $}
Предположим, что существует $2^n$ функций $h_s \in
C^1(\mathbb{R}\times \mathbb{R}^n \times [0,1],\mathbb{R}^n),
\quad s \in \{-1,1 \}^n$ таких, что\\ \centerline{$
h(t,\xi,\varepsilon)=h_s(t,\xi,\varepsilon), \quad
(t,\xi,\varepsilon)\in \mathbb{R}\times \mathbb{R}^n \times [0,1],
\quad s \in \{-1,1 \}^n $} Предположим, что $D \in
C^0(\mathbb{R}\times \mathbb{R}^n \times [0,1],\mathbb{R}^n)$ и
функция $D_1(\cdot,v_0,0)\cdot...\cdot D_n(\cdot,v_0,0)$ имеет
конечное число нулей на $[0,T]$. Тогда функция $
  (t,v,\varepsilon) \mapsto h(t,D(t,v,\varepsilon),\varepsilon)
$
удовлетворяет требованию интегральной дифференцируемости.
\end{lem}
\noindent {\bf Доказательство.} Пусть $0\leqslant t_1 \leqslant
t_2 \leqslant \dots \leqslant t_k \leqslant T$ -- все нули функции
$D_1(\cdot,v_0,0)\cdot...\cdot D_n(\cdot,v_0,0)$ на отрезке
$[0,T]$. Зафиксируем $\gamma
>0$ и положим
\centerline{$
  M=\left(\bigcup_{i=1}^k \left(t_i-\frac{\gamma}{2k},t_i+\frac{\gamma}{2k}\right)\right)\cap [0,T]
$} то есть
\begin{equation}\label{PROT}
D_1(\cdot,v_0,0)\cdot...\cdot D_n(\cdot,v_0,0)\not=0\quad\mbox{для
всех} \ t\in[0,T]\backslash M \end{equation} Покажем, что
существует $\delta>0$ такое, что
$D_1(t,v,\varepsilon)\cdot...\cdot D_n(t,v,\varepsilon)\not=0$ для
всех $t \in [0,T]\setminus M, \quad \|v-v_0\|<\delta, \quad
\varepsilon \in [0,\delta].$ Действительно, предположив противное,
получаем последовательность $(t_m,v_m,\varepsilon_m) \to
(t_0,v_0,0)$ при $m \to \infty$, где $t_0 \in [0,T]\setminus M$ и
такую,  что $D_1(t_m,v_m,\varepsilon_m)\cdot \dots \cdot
D_n(t_m,v_m,\varepsilon_m)=0$ при $m \in \mathbb{N}$, в чем
противоречие с (\ref{PROT}). Таким образом, при каждом $t \in
[0,T]\setminus M, \quad \|v-v_0\|<\delta$ и $\varepsilon
\in[0,\delta]$ найдется $s_{t,v,\varepsilon}\in \{-1,1\}^n$ такое,
что $D(t,v,\varepsilon) \in R^n_{s_{t,v,\varepsilon}}$ и, значит,
функция $h(t,D(t,\cdot,\varepsilon),\varepsilon)$ дифференцируема
в точке $v$. Заметим, что $\delta>0$ может быть уменьшено еще и
настолько, что \centerline{$ s_{t,v,\varepsilon}=s_{t,v_0,0} \quad
\mbox{при всех} \quad t \in [0,T]\setminus M, \quad
\|v-v_0\|<\delta, \quad \varepsilon \in [0,\delta] $} Аналогично,
доказывая от противного и выделяя соответствующую сходящуюся
подпоследовательность, приходим к заключению, что $ D(t_0,v_0,0)
\in R^n_{\lim\limits_{m \to \infty}s_{t_m,v_m,\varepsilon_m}},$
где $\lim_{m \to \infty}s_{t_m,v_m,\varepsilon_m} \neq
s_{t_0,v_0,0}$ и $t_0 \in [0,T] \setminus M.$ Но, по определению,
$R^n_{s_1} \cap R^n_{s_2}=\emptyset$ при $s_1 \neq s_2$.
Полученное противоречие завершает доказательство.\qed

Частный случай функций $D_1$ и $D_2$ при $n=2$ изучался ранее в
\cite{blm} и в работах А.~Буйка.

Оказывается, для систем (\ref{psb}), в которых правая часть не
везде дифференцируема, но липшицева и удовлетворяет требованию
интегральной дифференцируемости, утверждение второй теоремы
Н.~Н.~Боголюбова \cite[Теорема~II]{bog} полностью справедливо.
Соответствующий результат установлен автором в \cite{blm}, но
некоторые похожие идеи неявно использовались еще в классических
работах по динамике подвесных мостов (см. напр. \cite{glover}).
Более точно, имеет место теорема.
\begin{thm}\label{thm2}{\bf \hskip-0.2cm.}
Пусть функция $h \in C^0(\mathbb{R}\times \mathbb{R}^n \times
[0,1],\mathbb{R}^n)$ и $D \in C^0(\mathbb{R}\times \mathbb{R}^n
\times [0,1],\mathbb{R}^n)$ $T$-периодичны по времени и
удовлетворяют условиям леммы~\ref{lem2}. Рассмотрим систему
\begin{equation}\label{psb}
\dot{A}=\varepsilon h(t,D(t,A,\varepsilon),\varepsilon)
\end{equation}
и соответствующую функцию усреднения
\begin{equation*}
\overline{h}({A})=\frac{1}{T}\int^T_0 h(\tau,D(\tau,A,0),0)d\tau
\end{equation*}
Пусть $A_0 \in \mathbb{R}^n$ -- некоторый нуль функции
$\overline{h}$, в окрестности некоторого $\overline{h}$ непрерывно
дифференцируема. Тогда, если все собственные значения матрицы
$\overline{h}'(A_0)$ отрицательны, то при всех достаточно малых
$\varepsilon
>0$ система (\ref{psb}) имеет единственное $T$- периодическое
решение $A_\varepsilon$ такое, что $A_\varepsilon (t) \to A_0$ при
$\varepsilon \to 0$. Более того, решение $A_\varepsilon$
асимптотически устойчиво.
\end{thm}
Доказательство следует из (\cite{blm}, теорема 2.5),
леммы~\ref{lem2} и замечания о том, что в условиях
леммы~\ref{lem2} функция $(t,v,\varepsilon) \mapsto
h(t,D(t,v,\varepsilon),\varepsilon)$ липшицева по $v$ равномерно
на $\mathbb{R}\times \mathbb{R}^n \times [0,1]$.

{\bf 4. Асимптотическая устойчивость и двухчастотность резонансных
колебаний двухмассного резонансного грохота.} Уравнения изучаемого
грохота (фиг.~1) выписываются (см. \cite{kru}, с.~133) как
\setcounter{subsection}{4} \setcounter{equation}{0}
\begin{eqnarray}
&&\hskip-2cm m_1\ddot{x}_1+\varepsilon \dot{P}(x_1-x_2)(\dot{x}_1-\dot{x}_2)-\varepsilon k_0 (\dot{\eta}(t)-\dot{x}_1+\dot x_2)+ \nonumber\\
&&+P(x_1-x_2)-k_0(\eta(t)-x_1+x_2)=0 \label{GR1}\\
&&\hskip-2cm m_2\ddot{x}_2-\varepsilon \dot{P}(x_1-x_2)(\dot{x}_1-\dot{x}_2)+\varepsilon k_2\dot{x}_2 +\varepsilon k_0(\dot{\eta}(t)- \nonumber \\
&&-\dot{x}_1+\dot x_2)-P(x_1-x_2)+k_2x_2+k_0(\eta(t)-x_1+x_2)=0
\label{GR2}
\end{eqnarray}
Предположение о малости жесткости пружины приводит к функции $P$
вида $ P(\Delta)= \varepsilon \widehat{k}_1 \Delta^+, $ где $
\Delta^+=\max \{0,\Delta \}. $ Коэффициенты $k_0$ и $k_2$ удобно
представить в виде $ k_0=\overline{k}_0+\varepsilon
\widetilde{k}_0,$ $k_2=\overline{k}_2+\varepsilon
\widetilde{k}_2.$ Для достижения резонанса амплитуда, сообщаемая
электродвигателем, выбирается порядка $\varepsilon>0$, то есть
положим $ \overline{k}_0\eta(t)=\varepsilon \widetilde{\eta}(t). $
Для приведения полученной системы к стандартной форме принципа
усреднения используем результат раздела 2, функции которого
принимают в случае системы (\ref{GR1})-(\ref{GR2}) следующий вид:
$$
  M = \begin{pmatrix} a_{11} & a_{12} \\ a_{21} & a_{22}\end{pmatrix}=\begin{pmatrix} \overline{k}_0&-\overline{k}_0\\ -\overline{k}_0&\overline{k}_0+\overline{k}_2
  \end{pmatrix}, \ \begin{pmatrix} Q_1(\eps)\\
  Q_2(\eps)\end{pmatrix}=\eps\widehat{k}_1(x_1-x_2)\begin{pmatrix}1
  \\ -1\end{pmatrix}
$$
и $q_1(v)=q_2(v)=v^+.$ Функцию  (\ref{g}) запишем в удобном для
применения теоремы (\ref{thm2}) виде:
\begin{equation*}
g(t,A,\varepsilon)= \widetilde{g}(t,D(t,A,\varepsilon))
\end{equation*}
где
\begin{eqnarray*}
\widetilde{g}(t, \xi, \varepsilon)&=& \begin{pmatrix} -\frac{1}{m_1} (\varepsilon \widehat{k}_1 \xi_1)\\ -\frac{1}{m_2} (\varepsilon \widehat{k}_1 \xi_1)\\
-\widehat{k}_1(\xi_1)^+ +\widetilde{k}_0(-\xi_1)+\overline{k}_0(-\xi_2)+\widetilde{\eta}(t)+\varepsilon \widetilde{k}_0(-\xi_2)+\varepsilon k_0\dot{\eta}(t)\\
-\widehat{k}_1(\xi_1)^+  -\widetilde{k}_0(-\xi_1)-\overline{k}_0(-\xi_2)-\widetilde{\eta}(t)- \widetilde{k}_2\xi_4-\overline{k}_2(\xi_4-\xi_2)\\
-\varepsilon \widetilde{k}_0(-\xi_2)-\varepsilon k_0
\dot{\eta}(t)-\varepsilon \widetilde{k}_2(\xi_4-\xi_2)
\end{pmatrix}
\end{eqnarray*}
\begin{eqnarray*}
D(t,A,\varepsilon)&=& \begin{pmatrix} 1 & 0 & 0 & 0\\
-\eps^2\widehat{k}_1\frac{1}{m_1}+\eps^2\widehat{k}_1\frac{1}{m_2}
& 1 & 0 & 0\\
0 & 0 & 1 & 0\\
-\eps^2\widehat{k}_1\frac{1}{m_1} & 0 & 0 & 1
\end{pmatrix}\begin{pmatrix} \Omega_1(t)A-\Omega_2(t)A\\
\frac{1}{m_1}\Omega_3(t)A-\frac{1}{m_2}\Omega_4(t)A
\\ \Omega_1(t)A \\
\frac{1}{m_1}\Omega_3(t)A
\end{pmatrix}
\end{eqnarray*}
Так как дифференцируемость функции $(t,\xi,\varepsilon) \mapsto
\widetilde{g}(t,\xi,\varepsilon)$ нарушается только на
гиперплоскости $\xi_1=0$, то функция $
h(t,\xi,\varepsilon)=\Omega^{-1}(t)\widetilde{g}(t,\xi,\varepsilon)
$ удовлетворяет условиям леммы~\ref{lem2}. Так как строки матрицы
$\Omega(t)$ линейно независимы, то условиям леммы~\ref{lem2}
удовлетворяет также и функция $D.$ Для применения
теоремы~\ref{thm2} остается найти условия, при которых функция
\begin{equation}\label{FUFU}
t \mapsto
\Omega^{-1}(t)\widetilde{g}(t,D(t,A,\varepsilon),\varepsilon)
\end{equation}
периодична. Условие (\ref{A1}) приобретает вид $
-\overline{k}_0(-\overline{k}_0) <
-\overline{k}_0(-\overline{k}_2+\overline{k}_0), $ то есть
выполнено всегда. Хотя для периодичности функции (\ref{FUFU})
достаточно, чтобы частоты $\omega_1$ и $\omega_2$ из (\ref{Omega})
были соизмеримы, наибольший интерес представляет случай, когда $
\omega_2=l\omega_1$ для некоторого $l\in \mathbb{N}.$ Подставляя
даваемые формулой (\ref{LAM}) значения $\omega_1$ и $\omega_2,$
приходим к следующему условию периодичности функции (\ref{FUFU})
\begin{eqnarray}
&&l^2 \left(\frac{\overline{k}_2+\overline{k}_0}{m_2}+\frac{\overline{k}_0}{m_1}\right)^2=(1+l^2)^2 4\frac{\overline{k}_0 \overline{k}_2}{m_1m_2} \quad \mbox{для некоторого}  \quad l \in \mathbb{N} \label{CO1}\\
&&\widetilde{\eta} \left( t+\frac{2\pi}{\min \{\omega_1,\omega_2
\}}\right)=\widetilde{\eta}(t) \quad \mbox{при всех} \quad t \in
\mathbb{R} \label{CO2}
\end{eqnarray}
Алгоритм вычисления среднего от функции (\ref{FUFU}) определяется
числом $l$ в (\ref{CO1}) и функцией $\eta(t)$ в (\ref{CO2}). В
этой статье мы проводим это вычисление для случая $l=2$ и
$\eta(t)=\cos(\omega t),$ что позволяет заменить основное явление
-- двухчастотность колебаний грохота.

{\bf 4.1. Вычисление функции усреднения в случае $\omega_2=2w,
\omega_1=w.$ Доказательство наличия в порождающем решении
гармоники с частотой $2w.$}  Усредняя функцию (\ref{FUFU}) за
период $2\pi/\omega,$ при $\widehat{k}_1=0$
 получаем
\begin{equation*}
\overline{h}_0(A_{1C},A_{1S},A_{2C},A_{2S})= \begin{pmatrix}
-\alpha A_{1C}-\beta A_{1S}+\mu\\\beta A_{1C}-\alpha
A_{1S}\\-\gamma A_{2C}-\sigma A_{2S}\\ \sigma A_{2C}-\gamma A_{2S}
\end{pmatrix}
\end{equation*}
где
\begin{eqnarray*}
&&\alpha = \frac{\omega}{6} \left( \frac{\overline{k}_0 \overline{k}_2}{m_1m_2\omega^3}-\frac{\overline{k}_0}{m_1\omega}-\frac{\overline{k}_0+\overline{k}_2}{m_2\omega}+4\omega \right)\\
&&\beta =  \frac{\omega}{6} \left( \frac{4\widetilde{k}_0}{\overline{k}_0}+\frac{\overline{k}_0 \widetilde{k}_2}{m_1m_2\omega^4}-\frac{\widetilde{k}_0}{m_1\omega^2}-\frac{\widetilde{k}_0+\widetilde{k}_2}{m_2\omega^2} \right)\\
&&\gamma =\frac{\omega}{6} \left( \frac{4\overline{k}_0}{m_1\omega}+\frac{4\overline{k}_0}{m_2\omega}+\frac{4\overline{k}_2}{m_2\omega}-\frac{\overline{k}_0 \overline{k}_2}{m_1m_2\omega^3}-4\omega \right)\\
&&\sigma=\frac{\omega}{6} \left( \frac{2\widetilde{k}_0}{m_1\omega^2}+ \frac{2\widetilde{k}_0}{m_2\omega^2}+ \frac{2\widetilde{k}_2}{m_2\omega^2}-\frac{2\widetilde{k}_0}{\overline{k}_0}-\frac{\overline{k}_0 \widetilde{k}_2}{2m_1m_2\omega^4} \right)\\
&&\mu= \frac{1}{6m_1\omega^3} \left( \frac{1}{m_1}+\frac{1}{m_2}
\right) r - \frac{2}{3\overline{k}_0m_1\omega}r
\end{eqnarray*}
Условия (см. теорему~\ref{thm2}) существования у функции
$\overline{h}_0$ единственного нуля $A_0$ такого, что собственные
значения матрицы $(\overline{h}_0)'(A_0)$ отрицательны, приводит к
следующим предположениям:
\begin{equation}\label{LAGA}
\alpha>0, \gamma >0
\end{equation}
и соответствующий нуль дается формулой
\begin{equation}\label{AAAA}
(\widetilde{A}_{1C},\widetilde{A}_{1S},\widetilde{A}_{2C},\widetilde{A}_{2S})=
\left( \frac{\mu \alpha}{\alpha^2+\beta^2}, \frac{\mu
\beta}{\alpha^2+\beta^2},0,0 \right)
\end{equation}
Таким образом, приходим к следующему предложению.
\begin{pro}\label{pro1}{\bf \hskip-0.2cm.}
$(\widehat{k}_1=0)$ Пусть для $l=2$ выполнены предположения
(\ref{CO1})-(\ref{CO2}). Пусть выполнены условия (\ref{LAGA}).
Тогда, при $\widehat{k}_1=0$ и всех достаточно малых
$\varepsilon>0$, система (\ref{GR1})-(\ref{GR2}) имеет
единственное $\frac{2\pi}{\omega}$-периодическое решение
$(x_{1,\varepsilon},x_{2,\varepsilon})$ такое, что
\begin{equation*}
\begin{pmatrix} x_{1,\varepsilon}(t)\\x_{2,\varepsilon}(t) \end{pmatrix} \to
\begin{pmatrix} -\overline{k}_0 \frac{\mu \alpha}{\alpha^2+\beta^2}\sin (\omega t)+\overline{k}_0\frac{\mu \beta}{\alpha^2+\beta^2}\cos (\omega t)\\
(-\overline{k}_0+m_1\omega^2)\frac{\mu
\alpha}{\alpha^2+\beta^2}\sin (\omega
t)+(-\overline{k}_0+m_1\omega^2)\frac{\mu
\beta}{\alpha^2+\beta^2}\cos(\omega t) \end{pmatrix}
\end{equation*}
при $\varepsilon \to 0.$ Решение
$(x_{1,\varepsilon},x_{2,\varepsilon})$ асимптотически устойчиво.
\end{pro}
Пусть теперь $\widehat{k}_1 \ne 0$ и
$\overline{h}_{\widehat{k}_1}$ -- среднее значение соответствующей
функция (\ref{FUFU}). Прежде всего заметим, что существует
$\delta>0$ такое, что $\overline{h}_{\widehat{k}_1}$ непрерывно
дифференцируема при $\widehat{k}_1 \in [0,\delta]$ в
$\delta$-окрестности точки
$(\widetilde{A}_{1C},\widetilde{A}_{1S},\widetilde{A}_{2C},\widetilde{A}_{2S})$.
Утверждение следует из теоремы о неявной функции (\cite{kol},
Гл.~X, \S2) и того, что функция $t \mapsto
(\Omega_1(t)-\Omega_2(t))(\widetilde{A}_{1C},\widetilde{A}_{1S},\widetilde{A}_{2C},\widetilde{A}_{2S})^*$
пересекает $0$ трансверсально. Далее, вновь по теореме о неявной
функции, при всех достаточно малых ${\widehat{k}_1} >0$ функция
$\overline{h}_{\widehat{k}_1}$ имеет единственный нуль $(A_{1C,
\widehat{k}_1},A_{1S, \widehat{k}_1},A_{2C, \widehat{k}_1},A_{2S,
\widehat{k}_1})$ такой, что
\begin{equation}\label{CONV}
(A_{1C, \widehat{k}_1},A_{1S, \widehat{k}_1},A_{2C,
\widehat{k}_1},A_{2S, \widehat{k}_1}) \to \left( \frac{\mu
\alpha}{\alpha^2+\beta^2},\frac{\mu \beta}{\alpha^2+\beta^2},0,0
\right) \ \ \mbox{при} \ \ \widehat{k}_1 \to 0
\end{equation}
Из непрерывности производной функции
$\overline{h}_{\widehat{k}_1}$ следует, что при достаточно малых
${\widehat{k}_1} >0$ собственные значения матрицы
$(\overline{h}_{\widehat{k}_1})'(A_{1C, \widehat{k}_1},A_{1S,
\widehat{k}_1},A_{2C, \widehat{k}_1},A_{2S, \widehat{k}_1})$ имеют
отрицательные вещественные части. Уменьшим $\delta>0$ еще и так,
что
\begin{equation}\label{TCHTO}
A_{1C, \widehat{k}_1}>0 \quad \mbox{и} \quad A_{1S,
\widehat{k}_1}>0 \quad \mbox{при} \quad \widehat{k}_1 \in
[0,\delta]
\end{equation}
и покажем, что
\begin{equation}\label{SPRA}
A_{2C, \widehat{k}_1} \neq 0 \quad \mbox{или} \quad A_{2S,
\widehat{k}_1} \neq 0 \quad \mbox{при любом} \quad \widehat{k}_1
\in [0,\delta]
\end{equation}
Предположим, что это не так, то есть $ A_{2C, \widehat{k}_1}
=A_{2S, \widehat{k}_1} = 0$ при некотором $\widehat{k}_1 \in
[0,\delta].$ Для вычисления функции усреднения
$\overline{h}_{\widehat{k}_1}$ в точке
$\widetilde{\widetilde{A}}=(A_{1C, \widehat{k}_1},A_{1S,
\widehat{k}_1},0,0)$ вычислим
$(\Omega_1(t)\widetilde{\widetilde{A}}-\Omega_2(t)\widetilde{\widetilde{A}})^+$.
Имеем
\begin{equation}\label{Omega2}
\Omega_1(t)\widetilde{\widetilde{A}}-\Omega_2(t)\widetilde{\widetilde{A}}=(a_{12}+a_{11}-m_1
\omega^2)(A_{1C, \widehat{k}_1}\sin(\omega t)+A_{1S,
\widehat{k}_1}\cos(\omega t)).
\end{equation}
Два последовательных нуля функции (\ref{Omega2}) определяются
формулами
\begin{equation*}
t_1=-t_2,\quad t_2=\frac{1}{\omega} \arccos \left( \frac{A_{1C,
\widehat{k}_1}}{\sqrt{(A_{1C, \widehat{k}_1})^2+(A_{1S,
\widehat{k}_1})^2}} \right)
\end{equation*}
при этом
$\left(\Omega_1(\cdot)\widetilde{\widetilde{A}}-\Omega_2(\cdot)\widetilde{\widetilde{A}}\right)'(t_1)=-m_1\omega^3
\sqrt{(A_{1C, \widehat{k}_1})^2+(A_{1S, \widehat{k}_1})^2}, $
откуда следует, что
$(\Omega_1(t)\widetilde{\widetilde{A}}-\Omega_2(t)\widetilde{\widetilde{A}})^+=
\Omega_1(t)\widetilde{\widetilde{A}}-\Omega_2(t)\widetilde{\widetilde{A}},$
при $t \in [t_1,t_2],$ и
$(\Omega_1(t)\widetilde{\widetilde{A}}-\Omega_2(t)\widetilde{\widetilde{A}})^+=0,$
при $t\in[t_2,t_1+2\pi/\omega].$ Полученный вид функции $t \mapsto
(\Omega_1(t)\widetilde{\widetilde{A}}-\Omega_2(t)\widetilde{\widetilde{A}})^+$
позволяет вычислить функцию $\overline{h}_{\widehat{k}_1}(A_{1C,
\widehat{k}_1},A_{1S, \widehat{k}_1},0,0)$ и прийти к заключению,
что последние две строки системы
$\overline{h}_{\widehat{k}_1}(A_{1C, \widehat{k}_1},A_{1S,
\widehat{k}_1},0,0)=0$ записываются в виде
\begin{eqnarray}
& &\widehat{k}_1
\frac{-\overline{k}_0(m_1+m_2)+m_1m_2\omega^2}{18\overline{k}_0m_1m_2\omega^2}
\begin{pmatrix}A_{1S, \widehat{k}_1} & A_{1S, \widehat{k}_1} \\
-A_{2C, \widehat{k}_1} & A_{2C, \widehat{k}_1}\end{pmatrix}
 \left( \frac{3A_{1S, \widehat{k}_1}}{\sqrt{(A_{1C, \widehat{k}_1})^2+(A_{1S, \widehat{k}_1})^2}}, \right. \nonumber \\
& & \left.  \sin \left( 3\arccos \left( \frac{3A_{1C,
\widehat{k}_1}}{\sqrt{(A_{1C, \widehat{k}_1})^2+(A_{1S,
\widehat{k}_1})^2}} \right) \right) \right)^*=0 \label{U1}
\end{eqnarray}
Если
\begin{equation}\label{NUS}
m_1m_2\omega^2 \neq \overline{k}_0(m_1+m_2)
\end{equation}
то из соотношений (\ref{U1}) получаем, что по крайней мере одно из
чисел $A_{1S, \widehat{k}_1}$ и $A_{1C, \widehat{k}_1}$ равно
нулю, что противоречит (\ref{TCHTO}) и доказывает справедливость
(\ref{SPRA}). Таким образом, приходим к следующему предложению.
\begin{pro}\label{pro2}{\bf \hskip-0.2cm.}
$(\widehat{k}_1 \neq 0)$. Пусть выполнены условия
предложения~\ref{pro1} и условие (\ref{NUS}). Тогда существует
$\delta>0$, при котором каждому $\widehat{k}_1 \in (0,\delta]$
соответствует $\varepsilon_0>0$ такое, что при всех $\varepsilon
\in (0,\varepsilon_0]$ система (\ref{GR1})-(\ref{GR2}) имеет
единственное $T$-периодическое решение $(x_{1,\varepsilon},
x_{2,\varepsilon})$ такое, что
\begin{equation*}
\begin{pmatrix} x_{1,\varepsilon}(t)\\ x_{2,\varepsilon}(t) \end{pmatrix} \hskip-0.2cm\to
\hskip-0.2cm\left(\begin{array}{l} \hskip-0.2cm-\overline{k}_0A_{1C, \widehat{k}_1}\sin(\omega t)- \overline{k}_0A_{1S, \widehat{k}_1}\cos(\omega t)-\\
\hskip-0.2cm-\overline{k}_0A_{2C, \widehat{k}_1}\sin(2\omega t)-\overline{k}_0A_{2S, \widehat{k}_1} \cos(2\omega t)\\
\hskip-0.2cm(m_1\omega^2-\overline{k}_0)A_{1C, \widehat{k}_1}\sin(\omega t)+(m_1\omega^2-\overline{k}_0)A_{1S, \widehat{k}_1}+\\
\hskip-0.2cm+(4m_1\omega^2-\overline{k}_0)A_{2C,
\widehat{k}_1}\sin(2\omega t)+(4m_1\omega^2-\overline{k}_0)A_{2S,
\widehat{k}_1}\cos(2\omega t)
\end{array}\hskip-0.2cm\right)\hskip-0.2cm
\end{equation*}
при $\varepsilon \to 0$, где либо $A_{2C, \widehat{k}_1} \neq 0$,
либо $A_{2S, \widehat{k}_1} \neq 0$. При этом $(A_{1C,
\widehat{k}_1},A_{1S, \widehat{k}_1},A_{2C, \widehat{k}_1},A_{2S,
\widehat{k}_1}) \to (\frac{\mu \alpha}{\alpha^2+\beta^2},\frac{\mu
\beta}{\alpha^2+\beta^2},0,0)$ при $\widehat{k}_1 \to 0$ и решение
$(x_{1,\varepsilon},x_{2,\varepsilon})$ асимптотически устойчиво.
\end{pro}

{\bf 4.2. Пример конкретных параметров резонансного грохота,
приводящих к двухчастотным колебаниям.}  В этом подразделе
приводится пример конкретных параметров грохота, удовлетворяющих
всем условиям предложения (\ref{pro2}). Из (\ref{CO1}) имеем $
\overline{k}_2=\overline{k}_0 \cdot ({m_2}/{m_1}) \cdot
({1}/{4l^2}) \left( (1+l^2)- \sqrt{(1-l^2)-4l^2({m_1}/{m_2})}
\right)^2. $ В частности, приняв $m_1=11, m_2=64$, получаем $
\overline{k}_2= ({25}/{11})\overline{k}_0, $ на основании чего
зафиксируем $ \overline{k}_0=11, \overline{k}_2=25. $ Подставляя
полученные значения в формулу (\ref{LAM}), получаем $
\omega=\omega_1={\sqrt{5}}/{4},$
$2\omega=\omega_2={\sqrt{5}}/{2}.$ Соответствующие значения
$\alpha, \beta, \gamma, \sigma, \eta$ находятся как
$\alpha={\sqrt{5}\pi}/{4},$
$\beta={\pi}\widetilde{k}_0/132+{11\pi}\widetilde{k}_2/300,$
$\gamma=\sqrt{5}\pi,$
$\sigma={\pi}\widetilde{k}_0/{6}+{\pi}\widetilde{k}_2/{150},$
$\eta=-{4\pi}r/{1815},$ в частности, условие (\ref{LAGA})
выполнено. По формуле (\ref{AAAA}) находим
$(\widetilde{A}_{1C},\widetilde{A}_{1S})=-r\left(3403125+625(\widetilde{k}_0)^2+6050\widetilde{k}_0\widetilde{k}_2+14641(\widetilde{k}_2)^2\right)^{-1}(6000\sqrt{5},80(25\widetilde{k}_0+121\widetilde{k}_2)/11).$
Наконец, условие (\ref{NUS}), очевидно, выполнено и предложение
(\ref{pro2}) позволяет сформулировать следующий результат.

\begin{figure}[h]
\begin{center}
  \includegraphics[scale=1]{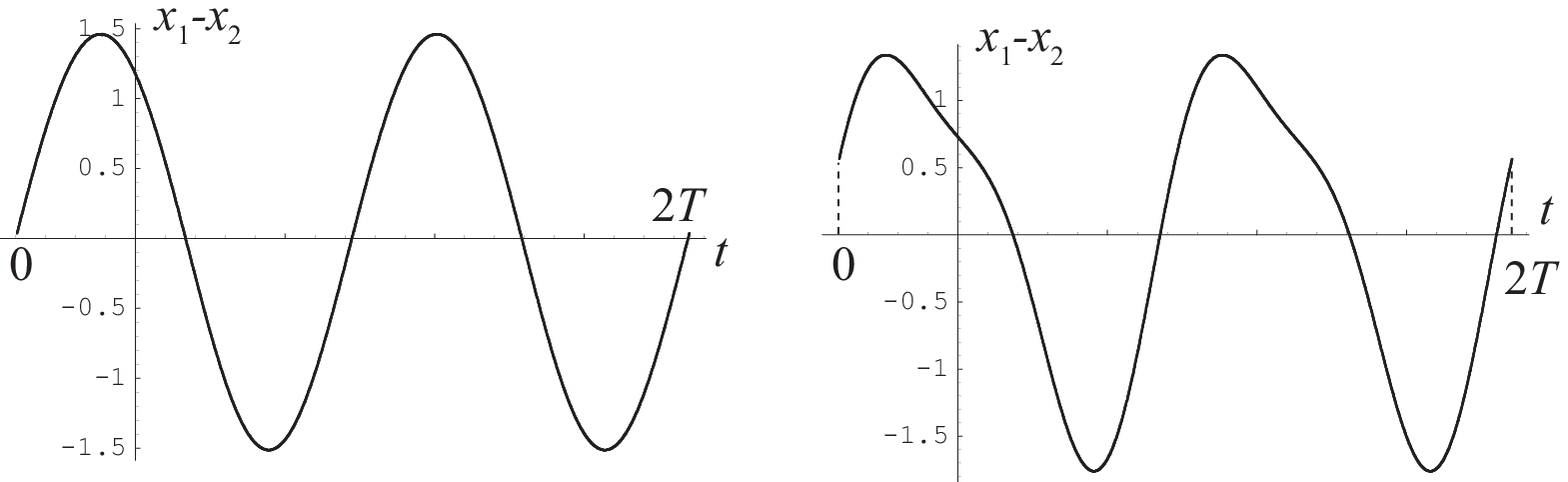}\\
  \caption{\footnotesize Численное моделирование асимптотически устойчивых решений системы (\ref{GR1})-(\ref{GR2}) при
  $m_1=11,$ $m_2=64,$ $\overline{k}_0=11,$ $\overline{k}_2=25,$ $\widetilde{k}_0=0,$ $\widetilde{k}_2=0,$
  $\eps=0.001,$ $r=10.$ Слева: $\widehat{k}_1=0,$ справа: $\widehat{k}_1=25.$}
  \end{center}\vskip-0.5cm
\end{figure}

\begin{pro}{\bf \hskip-0.2cm.}  Пусть $m_1=11,$ $m_2=64,$ $\overline{k}_0=11,$
$\overline{k}_2=25$ и $\widetilde{k}_0,\widetilde{k}_2,r>0$ --
произвольны. Тогда существует $\delta>0,$ при котором каждому
$\widehat{k}_1\in(0,\delta]$ соответствует $\eps_0>0$ такое, что
при всех $\eps\in(0,\eps_0]$ система (\ref{GR1})-(\ref{GR2}) имеет
единственное $8\pi/\sqrt{5}$-периодическое решение
$(x_{1,\eps},x_{2,\eps})$ такое, что
$$
\left(\begin{array}{l} \hskip-0.2cm x_{1,\eps}(t)\\
\hskip-0.2cm
x_{2,\eps}(t)\end{array}\hskip-0.2cm\right)\hskip-0.2cm\to\hskip-0.2cm
\left(\begin{array}{r} -11 A_{1C,\widehat{k}_1}
\sin\left({\sqrt{5}t}/{4}\right)-11 A_{1S,\widehat{k}_1}
\cos\left({\sqrt{5}t}/{4}\right)-\hskip-0.2cm\\-11
A_{2C,\widehat{k}_1} \sin\left({\sqrt{5}t}/{2}\right)-11
A_{2S,\widehat{k}_1} \cos\left({\sqrt{5}t}/{2}\right)\hskip-0.2cm\\
\hskip-0.2cm-({121}/{16}) A_{1C,\widehat{k}_1}
\sin\left({\sqrt{5}t}/{4}\right)-({121}/{16}) A_{1S,\widehat{k}_1}
\cos\left({\sqrt{5}t}/{4}\right)+\hskip-0.2cm\\+({11}/{4})
A_{2C,\widehat{k}_1} \sin\left({\sqrt{5}t}/{2}\right) +({11}/{4})
A_{2S,\widehat{k}_1}
\cos\left({\sqrt{5}t}/{2}\right)\hskip-0.2cm\end{array}\right)\hskip-0.2cm
$$
при $\eps\to 0,$ где либо $A_{2C,\widehat{k}_1}\not=0,$ либо
$A_{2S,\widehat{k}_1}\not=0.$ При этом
$\left(A_{1C,\widehat{k}_1},A_{1S,\widehat{k}_1},A_{2C,\widehat{k}_1},A_{2S,\widehat{k}_1}\right)\to
\left(\widetilde{A}_{1C},\widetilde{A}_{1S},0,0\right)$ при
$k_1\to 0$  и решение $(x_{1,\eps},x_{2,\eps})$ асимптотически
устойчиво.
\end{pro}

Представленное на фиг.~2 численное моделирование подтверждает, что
при $\widehat{k}_1\not=0$ колебания системы
(\ref{GR1})-(\ref{GR2}) действительно имеют неодночастотный
характер.

Работа поддержана грантом BF6M10 Роснауки и CRDF (программа BRHE)
и грантом MK-1620.2008.1 Президента РФ молодым кандидатам наук.
Исследования проведены в ходе стажировки автора в Институте
Проблем Управления РАН под руководством проф.~В.~Н.~Тхая и
финансируемой грантом РФФИ 08-01-90704-моб\_ст. Исследование из
раздела 3 поддержано грантом 09-01-90407-Укр\_ф\_а РФФИ, остальное
исследование -- грантом 09-01-00468-а. Автор благодарит
И.~С.~Мартынову за помощь в компьютерном наборе статьи.



\def\refname{\center \vskip-1cm ЛИТЕРАТУРА}

Воронеж \hfill Поступила в редакцию

e-mail: omakarenkov@math.vsu.ru \hfill 30.10.2008

\end{document}